\documentclass[12pt,a4paper]{article}
\usepackage{mathrsfs}
\usepackage{amssymb}
\usepackage{amsmath}

\newtheorem{theorem}{Theorem}[section]
\newtheorem{definition}[theorem]{Definition}
\newtheorem{lemma}[theorem]{Lemma}
\newtheorem{corollary}[theorem]{Corollary}

\date{}

\begin{document}

\title{Some remarks for the Akivis algebras and the Pre-Lie algebras\footnote{Supported by the
NNSF of China (Nos. 10771077; 10911120389).}}
\author{Yuqun Chen and Yu Li\\
{\small \ School of Mathematical Sciences, South China Normal
University}\\
{\small Guangzhou 510631, P. R. China}\\
{\small Email: yqchen@scnu.edu.cn}\\
{\small LiYu820615@126.com}}

\maketitle \noindent\textbf{Abstract:} In this paper, by using the
Composition-Diamond lemma for non-associative algebras invented by
A. I. Shirshov in 1962, we give Gr\"{o}bner-Shirshov bases for  free
Pre-Lie algebras and the universal enveloping non-associative
algebra of an Akivis algebra, respectively. As applications, we show
I.P. Shestakov's result that any Akivis algebra is linear and D.
Segal's result that the set of all good words in $X^{**}$ forms a
linear basis of the free Pre-Lie algebra $PLie(X)$ generated by the
set $X$. For completeness, we give the details of the proof of
Shirshov's Composition-Diamond lemma for non-associative algebras.



\section{Introduction}

A.G. Kurosh \cite{k} initiated to study free non-associative
algebras over a field proving that any subalgebra of a free
non-associative algebra is free. His student, A.I. Zhukov, proved in
\cite{z} that the word problem is algorithmically decidable in the
class of non-associative algebras. Namely, he proved that word
problem is decidable for any finitely presented non-associative
algebra. A.I. Shirshov, also a student of Kurosh, proved in
\cite{S2,S09}, 1953, that any subalgebra of a free Lie algebra is
free. This theorem is now known as the Shirshov-Witt theorem (see,
for example, \cite{R}) for it was proved also by E. Witt  \cite{W}.
Some later, Shirshov \cite{S1,S09} gave a direct construction of a
free (anti-) commutative algebra and proved that any subalgebra of
such an algebra is again free (anti-) commutative algebra. Almost
ten years later, Shirshov came back to, we may say, the Kurosh
programme, and published two papers \cite{S3} and \cite{S4}. In the
former, he gave a conceptual proof that the word problem is
decidable in the class of (anti-) commutative non-associative
algebras. Namely, he created the theory that is now known as
Gr\"{o}bner-Shirshov bases theory for (anti-) commutative
non-associative algebras. In the latter, he did the same for Lie
algebras (explicitly) and associative algebras (implicitly). Their
main applications were the decidability of the word problem for any
one-relater Lie algebra, the Freiheitsatz (the Freeness theorem) for
Lie algebras, and the algorithm for decidability of the word problem
for any finitely presented homogeneous Lie algebra. The same
algorithm is valid for any finitely presented homogeneous
associative algebra as well. Shirshov's main technical discovery of
\cite{S4,S09} was the notion of composition of two Lie polynomials
and implicitly two associative polynomials. Based on it, he gave the
algorithm to construct a Gr\"{o}bner-Shirshov basis for any ideal of
a free Lie algebra. The same algorithm is valid in the associative
case. This algorithm is in general infinite  as well as, for
example, Knuth-Bendix algorithm \cite{knuth}. Shirshov proved that
if a Gr\"{o}bner-Shirshov basis of an ideal is recursive, then the
word problem for the quotient algebra is decidable. It follows from
Shirshov's Composition-Diamond lemma that it is valid for  free
non-associative, free (anti-) commutative, free Lie and free
associative algebras (see  \cite{S3, S4,S09}). Explicitly the
associative case was treated in the papers by L.A. Bokut \cite{b76}
and G. Bergman \cite{b}.

Independently,  B.  Buchberger in his thesis (1965) (see
\cite{bu70}) created the Gr\"obner bases theory  for the classical
case of commutative associative algebras. Also, H. Hironaka in his
famous paper \cite{H} did the same for (formal or convergent)
infinite series rather than polynomials. He called his bases as the
standard bases. This term is used until now as a synonym of
Gr\"obner (in commutative case) or Gr\"obner-Shirshov (in
non-associative and non-commutative cases) bases.

There are a lot of sources of the history of Gr\"obner and
Gr\"{o}bner-Shirshov bases theory (see, for example,  \cite{E,
 B1,B2,B3}).

In the present paper we are dealing with the Composition-Diamond
lemma for a free non-associative algebra, calling it as
non-associative Composition-Diamond lemma. Shirshov mentioned it in
 \cite{S3,S09} that all his results are valid for the case of free
non-associative algebras rather than free (anti-) commutative
algebras. For completeness, we prove this lemma  in Section 2 in
this paper. Then we apply this lemma to the universal enveloping
non-associative algebra of an Akivis algebra and Pre-Lie algebra to
obtain Gr\"{o}bner-Shirshov bases for such algebras, respecrtively.
In particular, as applications, we show I.P. Shestakov's result that
any Akivis algebra is linear (see \cite{Sh}) and D. Segal's result
that the set of all good words in $X^{**}$ forms a linear basis of
the free Pre-Lie algebra $PLie(X)$ generated by the set $X$ (see
\cite{Se94}).

An Akivis algebra is a vector space $V$ over a field $k$ endowed
with a skew-symmetric bilinear product $[x,y]$ and a trilinear
product (x,y,z) that satisfy the identity
$[[x,y],z]+[[y,z],x]+[[z,x],y]=(x,y,z)+(z,x,y)+(y,z,x)-(x,z,y)-(y,x,z)-(z,y,x)$.
These algebras were introduced in 1976 by M.A. Akivis \cite{A} as
tangent algebras of local analitic loops. For any (non-associative)
algebra $B$ one may obtain an Akivis algebra $Ak(B)$ by considering
in $B$ the usual commutator $[x,y]=xy-yx$ and associator
$(x,y,z)=(xy)z-x(yz)$. Let $\{e_i\}_I$ be a linear basis of an
Akivis algebra $A$. Then the nonassociative algebra $U(A)=M(
\{e_i\}_I| \ e_ie_j-e_je_i=[e_i,e_j], \
(e_ie_j)e_k-e_i(e_je_k)=(e_i,e_j,e_k), \ i,j,k \in I)$ given by the
generators and relations is the universal enveloping non-associative
algebra of $A$, where $[e_i,e_j]=\sum_{m}\alpha_{ij}^me_m,\
(e_i,e_j,e_k)=\sum_{n}\beta_{ijk}^ne_n \mbox{ and each }
 \alpha_{ij}^m,\beta_{ijk}^n\in k$. The linearity of $A$ means that
 $A$ is a subspace of $U(A)$ (see
\cite{Sh}). Remark also that any subalgebra of a free Akivis algebra
is again free (see \cite{SU}).

A Pre-Lie algebra $A$ over a field $k$ is a non-associative algebra
with identity:
$$
(x,y,z)=(x,z,y), \ \ x,y,z\in A.
$$

\section{Composition-Diamond lemma for non-associative algebras}

Let $X=\{x_i|i\in I \}$ be a set, $X^*$ the set of all associative
words $u$ in $X$, and $X^{\ast \ast }$ the set of all
non-associative words $(u)$ in $X$. Let $k$  be a field and $M(X)$
be a $k$-space spanned by $X^{**}$. We define the product of
non-associative words by the following way:
$$
(u)(v)=((u)(v)).
$$
Then $M(X)$ is a free non-associative algebra generated by $X$.

Let $I$ be a well-ordered set. We order $X^{**}$ by the induction on
the length $|((u)(v))|$ of the words $(u)$ and $(v)$ in $X^{**}$:
\begin{enumerate}
\item[(i)] \  If $|((u)(v))|=2$, then $(u)=x_i > (v)=x_j$ if and only if $i>j$.
\item[(ii)] \ If $|((u)(v))|>2$, then $(u)>(v)$ if and only if one of the
following cases holds:
\begin{enumerate}
\item[(a)] $|(u)|>|(v)|$.
\item[(b)] \ If $|(u)|=|(v)|$ and $(u)=((u_1)(u_2))$,
$(v)=((v_1)(v_2))$, then $(u_1)>(v_1)$ or ($(u_1)=(v_1)$  and
$(u_2)>(v_2)$).
\end{enumerate}
\end{enumerate}

It is easy to check that $>$ is a monomial ordering on $X^{**}$ in
the following sense:
\begin{enumerate}
\item[(a)] \ $>$ is a well ordering.
\item[(b)] \
$ (u)>(v)\Longrightarrow (u)(w)>(v)(w) \ and  \ (w)(u)>(w)(v) \
\mbox { for any } (w)\in X^{**}. $
\end{enumerate}
Such an ordering is called deg-lex (degree-lexicographical) ordering
and we use this ordering throughout this paper.

Given a polynomial $f\in M(X)$, it has the leading word $(\bar f)
\in X^{**}$ according to the deg-lex ordering on $X^{**}$ such that
$$
f=\alpha(\overline{f})+\sum{\alpha}_i(u_i),
$$
where $(\overline{f})>(u_i), \ \alpha , {\alpha}_i \in k, \ (u_i)\in
X^{**}$. We call $(\overline{f})$ the leading term of $f$. $f$ is
called monic if $\alpha=1$.

\ \

Let $S\subset M(X)$ be a set of monic polynomials, $s\in S$ and $(u)\in X^{**}$. We define $S$-word $%
(u)_s$ by induction:
\begin{enumerate}
\item[(i)] $(s)_s=s$ is an $S$-word of $S$-length
1.
\item[(ii)] If $(u)_s$ is an $S$-word of $S$-length k and $(v)$
is a non-associative word of length $l$, then
\begin{equation*}
(u)_s(v)\ and\ (v)(u)_s
\end{equation*}
are $S$-words of length $k+l$.
\end{enumerate}

Note that for any $S$-word $(u)_s=(asb)$, where $a,b\in X^*$, we
have $\overline{(asb)}=(a\bar{s}b)$.

Let $f,g$ be monic polynomials in $M(X)$. Suppose that there exist
$a,b\in {X^*} $ such that $(\bar f) =(a(\bar g) b)$. Then we define
the composition of inclusion
\begin{equation*}
(f,g)_{(\bar f)}=f-(agb).
\end{equation*}
It is clear that
\begin{equation*}
(f,g)_{(\bar f)}\in Id(f,g)\ \ and\ \ \overline{(f,g)_{(\bar
f)}}<(\bar f)
\end{equation*}
where $Id(f,g)$ is the ideal of $M(X)$ generated by $f,g$.

The composition $(f,g)_{(\bar f)}$ is trivial modulo $(S,(\bar f))$,
if
\begin{equation*}
(f,g)_{(\bar f)}=\sum\limits_i\alpha_i(a_is_ib_i)
\end{equation*}
where each $\alpha_i\in k,\ a_i,b_i\in X^*,\ s_i\in S,\ (a_is_ib_i)$
an $S$-word and $(a_i(\bar{s_i})b_i)<(\bar f)$. If this is the case,
then we write $(f,g)_{(\bar f)}\equiv 0\ mod (S,(\bar f))$. In
general, for $p,q\in M(X)$ and $(w)\in X^{**}$, we write
$$
p\equiv q\quad mod(S,(w))
$$
which means that $p-q=\sum\alpha_i (a_i s_i b_i) $, where each
$\alpha_i\in k,a_i,b_i\in X^{*},\ s_i\in S,\ (a_is_ib_i)$ an
$S$-word and $(a_i (\bar {s_i}) b_i)<(w)$.

\begin{definition} (\cite{S3,S09})
Let $S\subset M(X)$ be a nonempty set of monic polynomials and the
ordering $>$ defined as before. Then $S$ is called a
Gr\"{o}bner-Shirshov basis in $M(X)$ if any composition
$(f,g)_{(\bar f)}$ with $f,g\in S$ is trivial modulo $(S,(\bar f))$,
i.e., $(f,g)_{(\bar f)}\equiv0$ $mod(S,(\bar f))$.
\end{definition}

\begin{lemma}\label{3.7}
Let $(a_1s_1b_1),~(a_2s_2b_2)$ be $S$-words. If $S$ is a
Gr\"{o}bner-Shirshov basis in $M(X)$ and
$(w)=(a_1(\overline{s_1})b_1)=(a_2(\overline{s_2})b_2)$, then
\begin{equation*}
(a_1s_1b_1)\equiv (a_2s_2b_2)\ mod (S,(w)).
\end{equation*}

\end{lemma}

\noindent\textbf{Proof.} We have $a_1\bar s_1 b_1=a_2\bar s_2 b_2$
as associative words in the alphabet $X\cup \{\bar {s_1}, \bar
{s_2}\}$. There are two cases to consider.

Case 1. Suppose that subwords $\bar s_1$ and $\bar s_2$ of $w$ are
disjoint, say, $|a_2|\geq |a_1|+|\bar s_1|$. Then, we can assume
that
$$
a_2=a_1\bar s_1 c \ \ and  \ b_1=c\bar s_2 b_2
$$
for some $c\in X^*$, and so, $ w=(a_1(\bar s_1) c (\bar s_2) b_2). $
Now,
\begin{eqnarray*}
(a_1 s_1 b_1)-(a_2 s_2 b_2)&=&(a_1 s_1 c (\bar s_2)
b_2)-(a_1(\bar s_1 )c  s_2 b_2)\\
&=&(a_1 s_1 c ((\bar s_2) - s_2) b_2)+(a_1(s_1-(\bar s_1)) c  s_2
b_2).
\end{eqnarray*}
Since $(\overline{(\overline{s_2})-s_2})<(\bar s_2)$ and
$(\overline{s_1-(\overline{s_1})})<(\bar s_1)$, we conclude that
$$
(a_1 s_1 b_1)-(a_2 s_2 b_2)=\sum\limits_i
\alpha_i(u_is_1v_i)+\sum\limits_j \beta_j(u_js_2v_j)
$$
for some $\alpha_i,\beta_j\in k$, $S$-words $(u_is_1v_i)$ and
$(u_js_2v_j)$ such that $ (u_i(\bar s_1)v_i),(u_j(\bar s_2)v_j)<(w).
$ Thus,
$$
(a_1s_1b_1)\equiv (a_2s_2b_2)\ mod (S,(w)).
$$

Case 2. Suppose that the subword $\bar s_1$ of $w$ contains $\bar
s_2$ as a subword. We assume that
$$
(\bar s_1)=(a(\bar s_2)b), \ a_2=a_1a  \mbox{ and } b_2=bb_1, \mbox{
that is, } (w)=(a_1a(\bar s_2)bb_1)
$$
for some $S$-word $(a s_2 b)$. We have
\begin{eqnarray*}
(a_1 s_1 b_1)-(a_2 s_2 b_2)&=&(a_1 s_1 b_1)-(a_1 (a s_2 b) b_1)\\
&=&(a_1(s_1-(as_2b))b_1)\\
&=&(a_1(s_1,s_2)_{(\overline{s_1})}b_1).
\end{eqnarray*}
Since $S$ is a Gr\"{o}bner-Shirshov basis,
$(s_1,s_2)_{(\overline{s_1})}=\sum\limits_i\alpha_i(c_i s_i d_i)$
for some $\alpha_i\in k$, $S$-words $(c_is_id_i)$ with each
$(c_i(\bar s_i) d_i)<(\bar{s_1})$. Then,
\begin{eqnarray*}
&&(a_1 s_1 b_1)-(a_2 s_2 b_2)=(a_1(s_1,s_2)_{(\overline{s_1})}b_1)\\
&=&\sum\limits_i\alpha_i(a_1(c_is_id_i)b_1)=\sum\limits_j\beta_j(a_js_jb_j)
\end{eqnarray*}
for some $\beta_j\in k $, $S$-words $(a_js_jb_j)$ with each
$(a_j(\bar s_j )b_j)<(w)=(a_1(\bar {s_1})b_1)$. \\
Thus,
$$
(a_1s_1b_1)\equiv (a_2s_2b_2)\ mod (S,(w)).  \ \ \ \ \square
$$

\begin{lemma}\label{3.8}
Let $S\subset M(X)$ be a subset of monic polynomials and
$Irr(S)=\{(u)\in X^{**} |(u)\ne (a(\bar s) b),\ a,b\in X^*,\ s\in S
\mbox{ and } (as b) \mbox{ is an  } S\mbox{-word}\}$. Then for any
$f\in M(X)$,
\begin{equation*}
f=\sum\limits_{(u_i)\leq (\bar f) }\alpha_i(u_i)+
\sum\limits_{(a_j(\overline{s_j})b_j)\leq(\bar f)}\beta_j(a_js_jb_j)
\end{equation*}
where each $\alpha_i,\beta_j\in k, \ (u_i)\in Irr(S)$ and
$(a_js_jb_j)$ an $S$-word.
\end{lemma}
{\bf Proof.} Let  $f=\sum\limits_{i}\alpha_{i}(u_{i})\in{M(X)}$,
where $0\neq{\alpha_{i}\in{k}}$ and $(u_{1})>(u_{2})>\cdots$. If
$(u_1)\in{Irr(S)}$, then let $f_{1}=f-\alpha_{1}(u_1)$. If
$(u_1)\not\in{Irr(S)}$, then there exist some $s\in{S}$ and
$a_1,b_1\in{X^*}$, such that $(\bar f)=(u_1)=(a_1(\bar{s_1})b_1)$.
Let $f_1=f-\alpha_1(a_1s_1b_1)$. In both cases, we have
$(\bar{f_1})<(\bar{f})$. Then the result follows from the induction
on $(\bar{f})$. \ \ \ \ $\square$

\ \

The proof of the following theorem is analogous to one in Shirshov
\cite{S3}. For convenience, we give the details.

\begin{theorem}\label{t1}(A.I. Shirshov \cite{S3,S09}, Composition-Diamond lemma for non-associative algebras)
Let $S\subset M(X)$ be a nonempty set of monic polynomials, $Id(S)$
the ideal of $M(X)$ generated by $S$ and the ordering $>$ on
$X^{**}$ defined as before. Then the following statements are
equivalent.
\begin{enumerate}
\item [(i)] $S$ is a Gr\"{o}bner-Shirshov basis in $M(X)$.

\item [(ii)] $f\in Id(S)\Rightarrow (\bar f) =(a(\bar s )b)$ for some $s\in S\
and\ a,b\in X^*$, where $(as b)$ is an $S$-word.

\item [(iii)] $Irr(S)=\{(u)\in X^{**} |(u)\ne (a(\bar s) b)\ a,b\in X^*,\ s\in S \mbox{ and }
(as b) \mbox{ is an } S\mbox{-word}\}$ is a linear basis of the
algebra $M(X|S)$.
\end{enumerate}
\end{theorem}
{\bf Proof.} $(i)\Rightarrow (ii)$. \ Let $S$ be a
Gr\"{o}bner-Shirshov basis and $0\neq f\in Id(S)$. Then, we have
$$
f=\sum_{i=1}^n\alpha_i(a_is_ib_i)
$$
where each $\alpha_i\in k, \ a_i,b_i\in {X^*}, \ s_i\in S$ and $
(a_is_ib_i)$ an $S$-word. Let
$$
(w_i)=(a_i(\overline{s_i})b_i), \
(w_1)=(w_2)=\cdots=(w_l)>(w_{l+1})\geq\cdots
$$
We will use the induction on $l$ and $(w_1)$ to prove that
$(\overline{f})=(a(\overline{s})b)$ for some $s\in S \ \mbox{and} \
a,b\in {X^*}$.

If $l=1$, then
$(\overline{f})=\overline{(a_1s_1b_1)}=(a_1(\overline{s_1})b_1)$ and
hence the result holds. Assume that $l\geq 2$. Then, by Lemma
\ref{3.7}, we have
$$
(a_1s_1b_1)\equiv(a_2s_2b_2) \ \ mod(S,(w_1)).
$$
Thus, if $\alpha_1+\alpha_2\neq 0$ or $l>2$, then the result holds.
For the case $\alpha_1+\alpha_2= 0$ and $l=2$, we use the induction
on $(w_1)$. Now, the result follows.\\

$(ii)\Rightarrow(iii)$. Suppose that
$\sum\limits_{i}\alpha_i(u_i)=0$ in $M(X|S)$, where $\alpha_i\in k$,
$(u_i)\in {Irr(S)}$. It means that
$\sum\limits_{i}\alpha_i(u_i)\in{Id(S)}$. Then all $\alpha_i$ must
be equal to zero. Otherwise,
$\overline{\sum\limits_{i}\alpha_i(u_i)}=(u_j)\in{Irr(S)}$ for some
$j$ which contradicts (ii).

Now, by Lemma \ref{3.8}, (iii) follows.\\

$(iii)\Rightarrow(i)$. For any $f,g\in{S}$ , by Lemma \ref{3.8} and
(iii), we have $ (f,g)_{(\bar f)}\equiv0\ \ \ mod(S,(\bar f)). $
Therefore, $S$ is a Gr\"{o}bner-Shirshov basis. \ \ \ \ $\square$

\section{Gr\"{o}bner-Shirshov basis for universal enveloping algebra of an Akivis algebra}

In this section, we obtain a Gr\"{o}bner-Shirshov basis for
universal enveloping non-associative algebra of an Akivis algebra.

\begin{theorem}
Let $(A,+,[-,-],(-,-,-))$ be an Akivis algebra over a field $k$ with
a well-ordered $k$-basis $\{e_i|\ i\in I\}$. Let
$$
[e_i,e_j]=\sum\limits_{m}\alpha_{ij}^me_m,\
(e_i,e_j,e_k)=\sum\limits_{n}\beta_{ijk}^ne_n,
$$
where $\alpha_{ij}^m,\beta_{ijk}^n\in k$. We denote
$\sum\limits_{m}\alpha_{ij}^me_m$ and
$\sum\limits_{n}\beta_{ijk}^ne_n$ by $\{e_ie_j\}$ and
$\{e_ie_je_k\}$, respectively. Let
$$
U(A)=M( \{e_i\}_I| \ e_ie_j-e_je_i=\{e_ie_j\}, \
(e_ie_j)e_k-e_i(e_je_k)=\{e_ie_je_k\}, \ i,j,k \in I)
$$
be the universal enveloping non-associative algebra of $A$. Let
\begin{eqnarray*}
S&=&\{f_{ij}=e_ie_j-e_je_i-\{e_ie_j\} \ (i>j), \
g_{ijk}=(e_ie_j)e_k-e_i(e_je_k)-\{e_ie_je_k\} \ (i,j,k\in I),
\\
&&h_{ijk}=e_i(e_je_k)-e_j(e_ie_k)-\{e_ie_j\}
e_k-\{e_je_ie_k\}+\{e_ie_je_k\} \ (i>j, \ k\geq j) \}.
\end{eqnarray*}
Then \begin{enumerate}
\item[(i)] \ $S$ is a Gr\"{o}bner-Shirshov basis in $M(\{e_i\}_I)$.
\item[(ii)] \ $Irr(S)=\{u|u\in \{e_i|\ i\in I\}^{**}\ \mbox{and} \ \ \
u\ \mbox{does \ not \ contain \ one \ of \ the \ words}\ e_ie_j\
(i>j), (e_ie_j)e_k \ (i,j,k\in I), e_i(e_je_k)\ (i>j, \ k\geq j)
\mbox{ as \ a \ subword}\}$ is a linear basis of the universal
enveloping non-associative algebra  $U(A)$ of $A$.
\item[(iii)] \ $A$ can be embedded into its universal
enveloping non-associative algebra $U(A)$.
\end{enumerate}
\end{theorem}

\noindent\textbf{Proof.} (i). It is easy to check that
$$
\overline{f_{ij}}=e_ie_j \ (i>j), \ \overline{g_{ijk}}=(e_ie_j)e_k \
(i,j,k\in I), \ \overline{h_{ijk}}=e_i(e_je_k)\ (i>j, \ k\geq j).
$$
So, we have only two kinds of compositions to consider:
$$
(g_{ijk},f_{ij})_{(e_ie_j)e_k}  \ (i>j,j\leq k) \ \mbox{ and } \
(g_{ijk},f_{ij})_{(e_ie_j)e_k}  \ (i>j> k).
$$
For $(g_{ijk},f_{ij})_{(e_ie_j)e_k},  \ (i>j,j\leq k)$, we have,
$mod(S,(e_ie_j)e_k)$,
\begin{eqnarray*}
&&(g_{ijk},f_{ij})_{(e_ie_j)e_k} \\
&=&(e_je_i)e_k-e_i(e_je_k)+\{e_ie_j\}e_k-\{e_ie_je_k\} \\
&\equiv&-e_i(e_je_k)+e_j(e_ie_k)+\{e_ie_j\}e_k+\{e_je_ie_k\}-\{e_ie_je_k\}
\\
&\equiv&0.
\end{eqnarray*}
For $(g_{ijk},f_{ij})_{(e_ie_j)e_k},  \ (i>j> k)$, by noting that,
in $A$,
\begin{eqnarray*}
&&[[e_i,e_j],e_k]+[[e_j,e_k],e_i]+[[e_k,e_i],e_j]\\
&=&(e_i,e_j,e_k)+(e_k,e_i,e_j)
+(e_j,e_k,e_i)-(e_i,e_k,e_j)-(e_j,e_i,e_k)-(e_k,e_j,e_i),
\end{eqnarray*}
we have, $mod(S,(e_ie_j)e_k)$,
\begin{eqnarray*}
&&(g_{ijk},f_{ij})_{(e_ie_j)e_k}\\
&=&(e_je_i)e_k-e_i(e_je_k)+\{e_ie_j\}e_k-\{e_ie_je_k\} \\
&\equiv&-e_i(e_je_k)+e_j(e_ie_k)+\{e_ie_j\}e_k+\{e_je_ie_k\}-\{e_ie_je_k\}
\\
&\equiv&-e_i(e_ke_j)-e_i\{e_je_k\}+e_j(e_ie_k)+\{e_ie_j\}e_k+\{e_je_ie_k\}-\{e_ie_je_k\}
\\
&\equiv&e_j(e_ie_k)-e_k(e_ie_j)-\{e_ie_k\}e_j+\{e_ie_j\}e_k-e_i\{e_je_k\}
\\
&&-\{e_ke_ie_j\}+\{e_ie_ke_j\}+\{e_je_ie_k\}-\{e_ie_je_k\}
\\
&\equiv&e_j(e_ke_i)+e_j\{e_ie_k\}-e_k(e_je_i)-e_k\{e_ie_j\}-\{e_ie_k\}e_j+\{e_ie_j\}e_k
\\
&&-e_i\{e_je_k\}-\{e_ke_ie_j\}+\{e_ie_ke_j\}+\{e_je_ie_k\}-\{e_ie_je_k\}
\\
&\equiv&e_k(e_je_i)+\{e_je_k\}e_i+\{e_ke_je_i\}-\{e_je_ke_i\}
         +e_j\{e_ie_k\}-e_k(e_je_i)-e_k\{e_ie_j\}
\\
&&-\{e_ie_k\}e_j+\{e_ie_j\}e_k-e_i\{e_je_k\}-\{e_ke_ie_j\}+\{e_ie_ke_j\}+\{e_je_ie_k\}-\{e_ie_je_k\}
\\
&\equiv&
\{e_je_k\}e_i-e_i\{e_je_k\}+e_j\{e_ie_k\}-\{e_ie_k\}e_j+\{e_ie_j\}e_k
-e_k\{e_ie_j\}
\\
&&+\{e_ke_je_i\}+\{e_ie_ke_j\}+\{e_je_ie_k\}-\{e_je_ke_i\}-\{e_ke_ie_j\}-\{e_ie_je_k\}
\\
&\equiv&
\{e_je_k\}e_i-e_i\{e_je_k\}+\{e_ke_i\}e_j-e_j\{e_ke_i\}+\{e_ie_j\}e_k
-e_k\{e_ie_j\}
\\
&&+\{e_ke_je_i\}+\{e_ie_ke_j\}+\{e_je_ie_k\}-\{e_je_ke_i\}-\{e_ke_ie_j\}-\{e_ie_je_k\}
\\
&\equiv& \{\{e_je_k\}e_i\}+\{\{e_ke_i\}e_j\}+\{\{e_ie_j\}e_k\}
\\
&&+\{e_ke_je_i\}+\{e_ie_ke_j\}+\{e_je_ie_k\}-\{e_je_ke_i\}-\{e_ke_ie_j\}-\{e_ie_je_k\}
\\
&\equiv&0.
\end{eqnarray*}
Thus, $S$ is a Gr\"{o}bner-Shirshov basis  in $M(\{e_i\}_I)$.

(ii) follows from Theorem \ref{t1}.

(iii) follows directly from (ii).

This completes our proof. \ \ \ \ $\square$

\section{Gr\"{o}bner-Shirshov bases for free Pre-Lie algebras}

In this section, we represent the free Pre-Lie algebra by
considering the free non-associative algebra and give a
Gr\"{o}bner-Shirshov basis for a free Pre-Lie algebra. As a result,
we re-show that the set of all good words in $X^{**}$ forms a linear
basis of the free Pre-Lie algebra $PLie(X)$ generated by the set $X$
(see \cite{Se94}).

The proof of the following theorem is straightforward and we hence
omit the details.

\begin{theorem}\label{t4.1}
Let $M(X)$ be the free non-associative algebra generated by $X$ and
let
\begin{equation*}
S=\{((u)(v))(w)-(u)((v)(w))-((u)(w))(v)+(u)((w)(v)) \ | \
(u),(v),(w)\in X^{**} \ and \ (v)>(w)\}.
\end{equation*}
Then the algebra $M(X|S)=M(X)/Id(S)$ is the free Pre-Lie algebra
generated by $X$. \ \  $\square$
\end{theorem}

We now cite the definition of good words (see \cite{Se94}) in
$X^{**}$ by induction on length:
\begin{enumerate}
\item[1)]$x_i$ is a good word for any $x_i\in X$.

Suppose that we define good words of length $<n$.

\item[2)] A word $((v)(w))$ is called a good word if and
only if \begin{enumerate}

\item[(a)] both $(v)$ and $(w)$ are good words,

\item[(b)] if $(v)=((v_1)(v_2))$, then $(v_2)\leq(w)$.
\end{enumerate}
\end{enumerate}

We denote $(u)$ by $[u]$, if $(u)$ is a good word. Let
\begin{eqnarray*}
&&S_0=\{([u][v])[w]-[u]([v][w])-([u][w])[v]+[u]([w][v]) \ | \\
&& \ \ \ \ \ \ \ \ \ \ [u],[v],[w]\ are\ good\ words \ and \
[v]>[w]\}.
\end{eqnarray*}

\begin{lemma}\label{l4.2}
Let $W$ be the set consisting of all good words. Then
$$
Irr(S_0)=\{(u)\in X^{**} |(u)\ne (a(\bar s) b),\ a,b\in X^*,\ s\in
S_0 \mbox{ and } (as b) \mbox{ is an  } s\mbox{-word}\ \}=W.
$$
\end{lemma}

\noindent\textbf{Proof.} Suppose that $(u)\in Irr(S_0)$. We will
show that $(u)$ is a good word by using induction on $|(u)|=n$. If
$n=1$, then $ (u)=x_i $ which is already a good word. Let $n>1$ and
$(u)=((v)(w))$. This case has two subcases. By induction, we see
immediately  that $ (v),(w) $ are both good words.

Subcase 1. If $|(v)|=1$, then $(u)$ is a good word.

Subcase 2. If $|(v)|>1\ and \ (v)=((v_1)(v_2))$, then $(v_2)\leq
(w)$ for $(u)\in Irr(S_0)$. Hence $(u)$ is a good word.

It is clear that every good word is in $Irr(S_0)$ since every
subword of a good word is still  a good word. \ \ \ \ $\square$

\ \

The following lemma follows from Lemmas \ref{3.8} and \ref{l4.2}.
\begin{lemma}\label{l4.3}
In $M(X)$, any word $(u)$ has the following presentation:
$$
(u)=\sum\limits_{i}\alpha_i[u_i]+\sum\limits_{j}\beta_j(a_js_jb_j),
$$
where $\alpha_i, \beta_j\in k$, $[u_i]$ are goood words,
$(a_j{s_j}b_j)$ are $S_0$-words, $s_j\in S_0, \ [u_i],
{(a_j(\overline{s_j})b_j)\leq(u)}$. Moreover, each $[u_i]$ has the
same length as $(u)$. \ \ \ \ $\square$
\end{lemma}

\begin{lemma}\label{l4.4}
Suppose that $S$ and $S_0$ are the sets defined as above. Then in
$M(X)$, we have
$$
Id(S)=Id(S_0).
$$
\end{lemma}
\noindent\textbf{Proof.} Since $S_0$ is a subset of $S$, we only
need to prove that $M(X|S_0)$ is a Pre-Lie algebra. In fact, we only
need to prove that the following hold in $M(X|S_0)$,
\begin{equation*}
((u)(v))(w)-(u)((v)(w))-((u)(w))(v)+(u)((w)(v))=0
\end{equation*}
where $(u),(v),(w)\in X^{**}$ and $(v)>(w)$. By Lemma 4.3, it
suffices to prove that for any good words $[u],[v],[w]$ with
$[v]>[w]$,
\begin{equation*}
([u][v])[w]-[u]([v][w])-([u][w])[v]+[u]([w][v])=0.
\end{equation*}
This is trivial by the definition of $S_0$. \ \ \ \ $\square$

\begin{theorem}\label{t4.5}
Let the ordering $>$ be defined as before and
\begin{center}
 $S_0=\{([u][v])[w]-[u]([v][w])-([u][w])[v]+[u]([w][v]) \ |
~[v]>[w]\  and  \ [u],[v],[w]\ are\ good\ words\}.$
\end{center}
 Then $S_0$ is a Gr\"obner-Shirshov basis in
$M(X)$.
\end{theorem}

\noindent\textbf{Proof.} To simplify our notations, we use $u$ for
$[u]$ and $u_1u_2\cdots u_n$ for
$(((u_1u_2)\cdots)u_n)$.\\
Let
\begin{equation*}
f_{uvw}=uvw-u(vw)-uwv+u(wv)
\end{equation*}
where $u,v,w $ are good words and $v>w$. It is easy to check that
$\overline{f_{uvw}}=uvw$.

Suppose that $\overline{f_{u_1v_1w_1}}$ is a subword of
$\overline{f_{uvw}}$. Since $u,v,w$ are good words, we have
$u_1v_1w_1=uv,\ u=u_1v_1, \ v=w_1$ and  $v_1>w_1=v>w$. We will prove
that the composition $ (f_{uvw},f_{u_1v_1w_1})_{uvw} $ is trivial
modulo $(S_0, uvw)$.

Firstly, we  prove that the following statements hold
$mod(S_0,uvw)$:
\begin{enumerate}
\item[1)]$u_1(v_1w)v-u_1(v_1wv)-u_1v(v_1w)+u_1(v(v_1w))\equiv 0$,

\item[2)]$u_1wv_1v-u_1w(v_1v)-u_1wvv_1+u_1w(vv_1)\equiv 0$,

\item[3)]$u_1(wv_1)v-u_1(wv_1v)-u_1v(wv_1)+u_1(v(wv_1))\equiv 0$,

\item[4)]$u_1(v_1v)w-u_1(v_1vw)-u_1w(v_1v)+u_1(w(v_1v))\equiv 0$,

\item[5)]$u_1vv_1w-u_1v(v_1w)-u_1vwv_1+u_1v(wv_1)\equiv 0$,

\item[6)]$u_1(vv_1)w-u_1(vv_1w)-u_1w(vv_1)+u_1(w(vv_1))\equiv 0$,

\item[7)]$u_1(vw)v_1-u_1(vwv_1)-u_1v_1(vw)+u_1(v_1(vw))\equiv 0$,

\item[8)]$u_1v_1(wv)-u_1(v_1(wv))-u_1(wv)v_1+u_1(wvv_1)\equiv 0$.
\end{enumerate}

We only prove 1). 2)-8) can be similarly proved. Denote by
$g=u_1(v_1w)v-u_1(v_1wv)-u_1v(v_1w)+u_1(v(v_1w))$. By Lemma
\ref{l4.3}, we have
$$
v_1w=\sum\limits_{i}\alpha_iu_i+\sum\limits_{j}\beta_j(a_js_jb_j)
$$
where $u_i$ are good words, $(a_j{s_j}b_j)$ are $S_0$-words, $s_j\in
S_0, \ u_i, {(a_j(\overline{s_j})b_j)\leq v_1w}$. Moreover, each
$u_i$ has the same length as $v_1w$.

By noting that
$u_1(a_j\bar{s_j}b_j)v,u_1((a_j\bar{s_j}b_j)v),u_1v(a_j\bar{s_j}b_j),
u_1(v(a_j\bar{s_j}b_j))<uvw$, we have
$$
g\equiv\sum\limits_{i}\alpha_ig_i \ \ \ mod(S_0,uvw)
$$
where $g_i=u_1u_iv-u_1(u_iv)-u_1vu_i+u_1(vu_i)$. Now $g_i=0$ or
$\bar{g_i}<uvw$ implies that $g_i\equiv0 \  \ mod(S_0,uvw)$ and so
$g\equiv0 \  \ mod(S_0,uvw)$.

Secondly we have
\begin{align*}
(f_{uvw},f_{u_1v_1w_1})_{uvw}&=f_{uvw}-(f_{u_1v_1w_1})w\\
&=-u_1v_1(vw)-u_1v_1wv+u_1v_1(wv)+u_1(v_1v)w+u_1vv_1w-u_1(vv_1)w.
\end{align*}
Then by $1)$--$6)$, we have, $mod(S_0, uvw)$,
\begin{eqnarray*}
-u_1v_1wv&\equiv&-u_1(v_1w)v-u_1wv_1v+u_1(wv_1)v\\
&\equiv &-u_1((v_1w)v)-u_1v(v_1w)+u_1(v(v_1w))-u_1w(v_1v)-u_1wvv_1\\
&&+u_1w(vv_1)+u_1(wv_1v)+u_1v(wv_1)-u_1(v(wv_1))\\
&\equiv &-u_1((v_1w)v)-u_1v(v_1w)+u_1(v(v_1w))-u_1w(v_1v)-u_1wvv_1+u_1w(vv_1)\\
&&+u_1(w(v_1v))+u_1(wvv_1)-u_1(w(vv_1))+u_1v(wv_1)-u_1(v(wv_1)),\\
\\
u_1(v_1v)w&\equiv&u_1(v_1vw)+u_1w(v_1v)-u_1(w(v_1v))\\
&\equiv&u_1(v_1(vw))+u_1(v_1wv)-u_1(v_1(wv))+u_1w(v_1v)-u_1(w(v_1v)),\\
\\
u_1vv_1w &\equiv&u_1v(v_1w)+u_1vwv_1-u_1v(wv_1)\\
&\equiv &u_1v(v_1w)+(u_1(vw))v_1+u_1wvv_1-u_1(wv)v_1-u_1v(wv_1),\\
\\
-u_1(vv_1)w &\equiv&-u_1(vv_1w)-u_1w(vv_1)+u_1(w(vv_1))\\
&\equiv
&-u_1(v(v_1w))-u_1(vwv_1)+u_1(v(wv_1))-u_1w(vv_1)+u_1(w(vv_1)).
\end{eqnarray*}

So, by $7)$--$8)$, we have, $mod(S_0, uvw)$,
\begin{align*}
(f_{uvw},f_{u_1v_1w_1})_{uvw}\equiv&-u_1v_1(vw)+u_1(v_1(vw))+u_1(vw)v_1-u_1(vwv_1)\\
&+u_1v_1(wv)+u_1(wvv_1)-u_1(v_1(wv))-u_1(wv)v_1\\
\equiv& 0.
\end{align*}

This completes the proof. \ \ \ \ $\square$

The following corollary follows from Lemmas \ref{l4.2}, \ref{l4.4}
and Theorems \ref{t4.1}, \ref{t4.5}, \ref{t1}.

\begin{corollary}
Let $M$ be the set of all good words in $X^{**}$ and $PLie(X)$ the
free Pre-Lie algebra over a field $k$ generated by $X$. Then the set
of all good words in $X^{**}$ is a linear basis of $PLie(X)$.
\end{corollary}

\ \

\noindent{\bf Acknowledgement}: The authors would like to express
their deepest gratitude to Professor L.A. Bokut for his kind
guidance, useful discussions and enthusiastic encouragement.

\end{document}